\documentclass[10pt,a4paper]{article}
\usepackage{amsmath,amssymb,amsthm,mathrsfs,dsfont}
\usepackage[margin=3cm]{geometry}
\usepackage[numbers,sort&compress]{natbib}
\usepackage{titlesec,hyperref}
\usepackage{color}
\usepackage{fancyhdr}
\usepackage{subfigure}
\usepackage[graphicx]{realboxes}
\usepackage{float}
\usepackage{empheq}
\titleformat{\subsection}
  {\it}
  {\thesubsection.\enspace}
  {1pt}
  {}

\newtheorem{theorem}{Theorem}[section]
\newtheorem{lemma}{Lemma}[section]
\newtheorem{remark}{Remark}[section]

\newtheorem{example}{Example}

\numberwithin{equation}{section}

\allowdisplaybreaks

\newcommand{\ep}{\mathrm{e}}

\begin{document}
\date{}
\title{Shock formation in traffic flow models with Arrhenius look-ahead interactions}
\author{
Jian Chen$^{1}$\thanks{chenj883@mail2.sysu.edu.cn (Jian Chen)},
\quad Zhaonan Luo$^{1}$\thanks{luozhn7@mail.sysu.edu.cn (Zhaonan Luo)},
\quad Zhaoyang Yin$^{1}$\thanks{mcsyzy@mail.sysu.edu.cn (Zhaoyang Yin)}
\\[2mm]
\small $^{1}$School of Science, Shenzhen Campus of Sun Yat-sen University,\\
\small Shenzhen 518107, China
}
\date{}
\maketitle

\begin{abstract}
We study finite-time shock formation in a nonlocal traffic flow model with Arrhenius look-ahead dynamics. For the constant interaction potential, we derive a time-dependent lower bound for the minimum slope and obtain a sufficient condition with logarithmic dependence on the magnitude of the minimum initial slope. For the linear interaction potential with look-ahead distance $\gamma$, we show that the condition $u'_0(x_0)>\frac{1}{\gamma}$ at a single point is sufficient for shock formation, which is independent of the minimum initial slope. Our results strictly improve the corresponding conditions of Lee and Liu \cite{r17}.

\par\smallskip

\end{abstract}
\noindent\rule{\textwidth}{0.5pt}

\section{Introduction}\label{sec:introduction}

In this paper, we consider the Cauchy problem for the nonlocal traffic flow model with Arrhenius look-ahead dynamics \cite{r24,r13}
\begin{align}\label{eq1}
\begin{cases}
u_t+\partial_x\bigl(u(1-u)\ep^{-\bar u}\bigr)=0,
& t>0,\quad x\in\mathbb R,\\
u(0,x)=u_0(x),&x\in\mathbb R,
\end{cases}
\end{align}
with either the constant interaction potential
\begin{align}\label{eq2}
\bar u(t,x)=\frac{1}{\gamma}\int_x^{x+\gamma}u(t,y)\,dy,
\end{align}
or the linear interaction potential
\begin{align}\label{eq3}
\bar u(t,x)=\frac{2}{\gamma}\int_x^{x+\gamma}
\left(1+\frac{x-y}{\gamma}\right)u(t,y)\,dy,
\end{align}
where $u(t,x)\in[0,1]$ denotes the normalized vehicle density and $\gamma>0$ is the look-ahead distance. Both interaction kernels are nonnegative and have unit mass. Eq.~\eqref{eq1} extends the Lighthill-Whitham-Richards model \cite{r21,r23} by incorporating the Arrhenius factor $\ep^{-\bar u}$, which accounts for the traffic ahead. The model with the constant interaction potential was derived from stochastic particle dynamics by Sopasakis and Katsoulakis \cite{r24}. Kurganov and Polizzi \cite{r13} studied numerical approximations and introduced the linear interaction potential \eqref{eq3}, which gives nearby vehicles greater influence on a car's speed than more distant vehicles. Sun and Tan \cite{r25} further proposed a cellular automaton model with two different types of look-ahead rules.

The mathematical analysis of related nonlocal traffic models includes well-posedness \cite{r1,r3,r7,r12}, numerical approximation \cite{r6},
stability \cite{r11}, and nonlocal-to-local limits \cite{r2,r4}. For the Arrhenius model \eqref{eq1}, a central question is whether the density gradient can become unbounded in finite time while the density remains bounded. For the constant potential, Li and Li \cite{r19} established local well-posedness and proved shock formation for several classes of initial data by comparing characteristic trajectories, see also \cite{r20}. Lee and Liu \cite{r17} subsequently derived explicit conditions for shock formation under the constant and linear interaction potentials by tracing the minimum and maximum slopes. Further developments include concavity changing fluxes \cite{r14}, look-ahead-behind interactions \cite{r15,r10}, and more general nonlocal conservation laws \cite{r16}. In particular, for the look-ahead-behind model with linear interaction potentials, Lee \cite{r15} obtained a condition for shock formation involving the initial slope at a single point.

For the infinite-range look-ahead interaction
\begin{align}\label{infinite}
\bar u(t,x)=\int_x^\infty u(t,y)\,dy, 
\end{align}
Lee and Tan \cite{r18} established a sharp threshold separating global regularity from finite-time shock formation. Hamori and Tan \cite{r8} extended this analysis to a more general class of fluxes. They also proved global regularity for a class of subcritical initial data in a second-order Aw-Rascle-Zhang model for the same look-ahead interaction \cite{r9}. The sharp threshold analyses in \cite{r18,r8} exploit the special structure $\bar u_x=-u$ of the nonlocal term \eqref{infinite}. For earlier work on critical threshold phenomena in nonlocal conservation laws, see Liu and Tadmor \cite{r22}.

In this paper, we establish refined sufficient conditions for finite-time shock formation in \eqref{eq1} with the constant and linear interaction potentials \eqref{eq2} and \eqref{eq3}. For the constant interaction potential, we obtain a time-dependent lower bound for the minimum slope and retain this time dependence in the differential inequality for the maximum slope. This yields a sufficient condition for shock formation with logarithmic dependence on the magnitude of the minimum initial slope. For the linear interaction potential, the identity for $\bar u_{xx}$ gives a differential inequality for $u_x$ along a single characteristic, which does not involve the minimum slope. The resulting condition for shock formation depends only on the initial slope at a single point. Our conditions strictly improve the results obtained by Lee and Liu \cite{r17}.

To state our results, we set (see Lemma \ref{L4})
\begin{align}\label{eq4}
m(t)=\underset{x\in\mathbb R}{\inf}u_x(t,x),
\qquad
M(t)=\underset{x\in\mathbb R}{\sup}u_x(t,x).
\end{align}
We first consider the constant interaction potential.

\begin{theorem}[Constant interaction potential]\label{thm:constant}
Let $u_0\in H^2(\mathbb R)$ satisfy $0\leq u_0\leq1$, and let
$T$ be the maximal existence time of the unique strong solution
of \eqref{eq1} with \eqref{eq2}. Set
$c_0=\max\{1,-\gamma m(0)\}$. If
\begin{align}\label{eq5}
M(0)>\frac{1}{\gamma}\left(\frac{1}{2}+\frac{\ep}{8}\log c_0\right),
\end{align}
then shock formation occurs in finite time with
\begin{align}\label{eq6}
0<T\leq\frac{2\ep\gamma}{3}
\log\left(1+\frac{3}{4\gamma\delta_0}\right),
\end{align}
where
\begin{align*}
\delta_0=M(0)-\frac{1}{\gamma}\left(\frac{1}{2}+\frac{\ep}{8}\log c_0\right)>0.
\end{align*}
\end{theorem}
\begin{remark}
\normalfont
\label{rem:constant-comparison}
\par\noindent (1) The sufficient condition for shock formation obtained in
\cite[Theorem~1.2]{r17} is
\begin{align}\label{eq7}
M(0)>\frac{1}{\gamma}\left(\frac{1}{2}+\frac{1}{4}\sqrt{2c_0+6}\right).
\end{align}
Since $\log c\leq\frac{2\sqrt c}{\ep}$ for $c\geq1$, we have
\begin{align*}
M(0)>\frac{1}{\gamma}
\left(\frac{1}{2}+\frac{\sqrt{2c_0+6}}{4}\right)
>\frac{1}{\gamma}
\left(\frac{1}{2}+\frac{\sqrt{c_0}}{4}\right)\geq\frac{1}{\gamma}
\left(\frac{1}{2}+\frac{\ep}{8}\log c_0\right).
\end{align*}
Thus \eqref{eq7} implies
\eqref{eq5}. Theorem~\ref{thm:constant} therefore improves the condition \eqref{eq7}.
\par\noindent (2) For fixed $\gamma$, the two conditions have different asymptotic behavior as $c_0\to\infty$. The threshold for $M(0)$ grows like $\sqrt{c_0}$ in \eqref{eq7} and only logarithmically with $c_0$ in \eqref{eq5}.
\end{remark}

\begin{example}
\normalfont
\label{ex:constant}
For example, we take
\begin{align*}
u_0(x)=\frac{3}{4}\sqrt{\frac{\ep}{2}}\,
\ep^{-\frac{x^2}{\gamma^2}}.
\end{align*}
Then $u_0\in H^2(\mathbb R)$, $0\leq u_0<1$,
$M(0)=\frac{3}{4\gamma}$ and $m(0)=-\frac{3}{4\gamma}$.
Thus $c_0=1$, and condition \eqref{eq5} is satisfied.
On the other hand, condition \eqref{eq7} becomes
\begin{align*}
M(0)>\frac{1}{\gamma}
\left(\frac{1}{2}+\frac{\sqrt2}{2}\right),
\end{align*}
which is not satisfied.
\end{example}

For the linear interaction potential, we obtain a criterion
involving only the initial slope at one point.

\begin{theorem}[Linear interaction potential]\label{thm:linear}
Let $u_0\in H^2(\mathbb R)$ satisfy $0\leq u_0\leq1$, and let
$T$ be the maximal existence time of the unique strong solution
of \eqref{eq1} with \eqref{eq3}. Assume that there is
$x_0\in\mathbb R$ such that
\begin{align}\label{eq8}
u'_0(x_0)>\frac{1}{\gamma}.
\end{align}
Then shock formation occurs in finite time with
\begin{align}\label{eq9}
0<T\leq\frac{\ep\gamma}{2}
\log\frac{\gamma u'_0(x_0)}{\gamma u'_0(x_0)-1}<\infty.
\end{align}
\end{theorem}

\begin{remark}
\normalfont
\label{rem:linear-comparison}
\par\noindent (1) Our criterion improves the shock formation result given in \cite[Theorem~1.3]{r17}. In fact, Lee and Liu obtained the following sufficient condition
\begin{align}\label{eq10}
M(0)>\frac{1}{\gamma}
\left(1+\frac{1}{2}\sqrt{6-\min\{-2,\gamma m(0)\}}\right).
\end{align}
It is obvious that
\begin{align*}
M(0)
&>\frac{1}{\gamma}
\left(1+\frac{1}{2}\sqrt{6-\min\{-2,\gamma m(0)\}}\right)\geq\frac{1+\sqrt{2}}{\gamma}
>\frac{1}{\gamma}.
\end{align*}
Therefore, Theorem~\ref{thm:linear} covers the condition \eqref{eq10}. Moreover, our condition does not depend on the minimum initial slope $m(0)$.
\par\noindent (2) For fixed $\gamma$, the constant $\frac{1}{\gamma}$ in the condition \eqref{eq8} could be optimized by a more delicate estimate in the proof of Theorem~\ref{thm:linear}. For simplicity, we omit the details in this result.
\end{remark}

\begin{example}
\normalfont
\label{ex:linear}
For example, we take
\begin{align*}
u_0(x)=\frac{1}{2}\ep^{-\frac{8\ep x^2}{\gamma^2}}.
\end{align*}
Then $u_0\in H^2(\mathbb R)$, $0\leq u_0\leq\frac{1}{2}$,
$M(0)=\frac{2}{\gamma}$ and $m(0)=-\frac{2}{\gamma}$.
Thus condition \eqref{eq8} is satisfied at
$x_0=-\frac{\gamma}{4\sqrt{\ep}}$.
On the other hand, condition \eqref{eq10} becomes
\begin{align*}
M(0)>\frac{1}{\gamma}(1+\sqrt2),
\end{align*}
which is not satisfied.
\end{example}

The remainder of the paper is organized as follows. In Section~\ref{sec:preliminaries}, we recall the local
well-posedness result and some useful lemmas. In Section~\ref{sec:proofs}, we prove our main results. 

\section{Preliminaries}\label{sec:preliminaries}
We begin with the local well-posedness result for \eqref{eq1}.
\begin{lemma}[\cite{r19,r17}]\label{lem:local}
Let $u_0\in H^2(\mathbb R)$.
Then, for either \eqref{eq2} or \eqref{eq3},
there exist a maximal existence time $T>0$ and a unique strong solution $u$ to \eqref{eq1} such that
\begin{align*}
u\in C([0,T);H^2(\mathbb R))
\cap C^1([0,T);H^1(\mathbb R)).
\end{align*}
Moreover, if $0\leq u_0\leq1$ and $T<\infty$, then
\begin{align}\label{eq11}
\underset{t\to T^-}{\limsup}
\|u_x(t,\cdot)\|_{L^\infty}=+\infty.
\end{align}
\end{lemma}
Consider the ordinary differential equation
\begin{align}\label{eq27}
q_t(t,x)=\bigl[(1-2u)\ep^{-\bar u}\bigr](t,q(t,x)),
\qquad q(0,x)=x.
\end{align}
Applying classical results in the theory of ordinary differential equations, \eqref{eq27} has a unique solution $q\in C^1([0,T)\times\mathbb R;\mathbb R)$. For $t\in[0,T)$, the map $q(t,\cdot)$ is an increasing $C^1$-diffeomorphism of $\mathbb R$ with
\begin{align*}
q_x(t,x)=\exp\left(\int_0^t
\partial_x\bigl[(1-2u)\ep^{-\bar u}\bigr]
(s,q(s,x))\,ds\right)>0.
\end{align*}

Next, the following lemma gives an a priori bound for the solution.
\begin{lemma}[\cite{r19,r17}]\label{lem:bounds}
Let $u_0\in H^2(\mathbb R)$ satisfy $0\leq u_0\leq1$,
and let $u$ be the solution given by Lemma~\ref{lem:local}.
For either \eqref{eq2} or \eqref{eq3}, we have
\begin{align*}
0\leq u(t,x)\leq1,
\qquad 0\leq\bar u(t,x)\leq1,
\qquad (t,x)\in[0,T)\times\mathbb R.
\end{align*}
Moreover
\begin{align}\label{eq12}
\ep^{-1}\leq\ep^{-\bar u(t,x)}\leq1.
\end{align}
\end{lemma}

Finally, we recall the following key lemma, which will be used in the proof of our main results.

\begin{lemma}[\cite{r5}]\label{L4}
Let $T>0$ and $u\in C^1([0,T];H^2(\mathbb R))$.
Then for every $t\in[0,T)$ there exists at least one point
$\xi(t)\in\mathbb R$ such that
\begin{align*}
m(t):=\underset{x\in\mathbb R}{\inf}\{u_x(t,x)\}=u_x(t,\xi(t)),
\end{align*}
and the function $m$ is almost everywhere differentiable on
$(0,T)$ with
\begin{align*}
\frac{dm(t)}{dt}=u_{tx}(t,\xi(t)),
\qquad\text{a.e. on }(0,T).
\end{align*}
\end{lemma}

\begin{remark}
\normalfont
The same statement holds for the supremum function
$M(t):=\underset{{x\in\mathbb R}}{\sup}\{u_x(t,x)\}$.
\end{remark}

\section{Proofs of the main results}\label{sec:proofs}
\subsection{Constant interaction potential}\label{subsec:constant}

\begin{proof}[\normalfont\bfseries Proof of Theorem~\ref{thm:constant}]
For the constant potential \eqref{eq2}, we have
\begin{align}\label{eq13}
\bar u_x=\frac{u(x+\gamma)-u(x)}{\gamma},
\qquad
\bar u_{xx}=\frac{u_x(x+\gamma)-u_x(x)}{\gamma}.
\end{align}
By Lemma~\ref{lem:bounds}, the first identity implies $-u\leq\gamma\bar u_x\leq1-u$. Differentiating \eqref{eq1} with respect to $x$, we obtain
\begin{align}\label{eq14}
u_{tx}+(1-2u)\ep^{-\bar u}u_{xx}
=\ep^{-\bar u}\bigl[2u_x^2+2(1-2u)\bar u_xu_x
+u(1-u)(\bar u_{xx}-\bar u_x^2)\bigr].
\end{align}
Let
\begin{align*}
m(t)=\inf\limits_{x\in\mathbb{R}}\{ u_x(t,x)\}=u_x\left(t,\xi_1(t) \right),\qquad M(t)=\sup\limits_{x\in\mathbb{R}}\{ u_x(t,x)\}=u_x\left(t,\xi_2(t) \right).
\end{align*}
By \eqref{eq13}, we have
\begin{align*}
\bar u_{xx}(t,\xi_1(t))\geq0,\qquad
\bar u_{xx}(t,\xi_2(t))\geq\frac{m-M}{\gamma}.
\end{align*}
Evaluating \eqref{eq14} at $\xi_1(t)$ and $\xi_2(t)$
and applying Lemma~\ref{L4}, we have for a.e. $t\in(0,T)$
\begin{align}
m'&\geq\ep^{-\bar u}
\bigl[2m^2+2(1-2u)\bar u_xm-u(1-u)\bar u_x^2\bigr],
\label{eq15}\\
M'&\geq\ep^{-\bar u}
\left[2M^2+2(1-2u)\bar u_xM
+u(1-u)\left(\frac{m-M}{\gamma}-\bar u_x^2\right)\right].
\label{eq16}
\end{align}
For $\gamma m\leq-1$, using $|\gamma\bar u_x|\leq1$, we obtain
\begin{align}\label{eqq1}
&\gamma^2\bigl[2m^2+2(1-2u)\bar u_xm
-u(1-u)\bar u_x^2\bigr]\nonumber\\
&\geq2\gamma^2m^2+2\gamma m|1-2u|-u(1-u)\nonumber\\
&\geq2\gamma^2m^2
+\gamma m\bigl[2|1-2u|+u(1-u)\bigr]\nonumber\\
&\geq2\gamma m(\gamma m+1),
\end{align}
where we have used $2|1-2u|+u(1-u)\leq2$ for $0\leq u\leq1$.
Combining \eqref{eqq1} with \eqref{eq15} and \eqref{eq12}, we obtain 
\begin{align}\label{eq17}
m'(t)\geq\frac{2}{\ep}m(t)\left(m(t)+\frac{1}{\gamma}\right)
\quad\text{for}\quad\gamma m(t)\leq-1.
\end{align}
To obtain a time-dependent lower bound for $m$, let $c$ solve
\begin{align}\label{eq18}
c'(t)=-\frac{2}{\ep\gamma}c(t)(c(t)-1),\qquad c(0)=c_0.
\end{align}
Then $1\leq c(t)\leq c_0$ for all $t\geq0$.
In view of $m(0)\geq-\frac{c_0}{\gamma}$, the comparison theorem for ordinary differential equations (see \cite[Lemma~16.4]{r26}) yields
\begin{align}\label{eq19}
m(t)\geq-\frac{c(t)}{\gamma},\qquad 0\leq t<T.
\end{align}
Substituting \eqref{eq19} into \eqref{eq16}, we obtain
\begin{align}\label{eq20}
\gamma^2M'\geq\ep^{-\bar u}\bigl\{
2\gamma^2M^2+2(1-2u)\gamma^2\bar u_xM
-u(1-u)\bigl[\gamma M+(\gamma\bar u_x)^2+c\bigr]\bigr\}.
\end{align}

We next estimate the right-hand side of \eqref{eq20} for $M\geq\frac{1}{2\gamma}$. By the symmetry $(u,\bar u_x)\mapsto(1-u,-\bar u_x)$, we may assume that $0\leq u\leq\frac{1}{2}$. Since the expression in braces in \eqref{eq20} is a concave function of $\bar u_x$, we only need to consider
\begin{align*}
\gamma\bar u_x=-u
\qquad\text{and}\qquad
\gamma\bar u_x=1-u.
\end{align*}
Subtracting the value at $\gamma\bar u_x=1-u$ from that at $\gamma\bar u_x=-u$, we obtain
\begin{align*}
&\Bigl(2\gamma^2M^2-2u(1-2u)\gamma M
-u(1-u)\bigl(\gamma M+u^2+c\bigr)\Bigr)\\
&\quad-\Bigl(2\gamma^2M^2+2(1-2u)(1-u)\gamma M
-u(1-u)\bigl(\gamma M+(1-u)^2+c\bigr)\Bigr)\\
&=(1-2u)\bigl(u(1-u)-2\gamma M\bigr)\\
&\leq-\frac{3}{4}(1-2u)\leq0.
\end{align*}
Thus it suffices to estimate the value at $\gamma\bar u_x=-u$, we have
\begin{align}\label{eq21}
&2\gamma^2M^2+2(1-2u)\gamma^2\bar u_xM
-u(1-u)\bigl[\gamma M+(\gamma\bar u_x)^2+c\bigr]\nonumber\\
&\geq2\gamma^2M^2+(5u^2-3u)\gamma M
-u^3(1-u)-cu(1-u)\nonumber\\
&\geq2\gamma^2M^2+(5u^2-3u)\gamma M
-\frac12u^2\gamma M-\frac14c\nonumber\\
&=2\gamma^2M^2+
\left[\frac92\left(u-\frac13\right)^2-\frac12\right]
\gamma M-\frac14c\nonumber\\
&\geq2\gamma^2M^2-\frac12\gamma M-\frac14c,
\end{align}
where we have used
\begin{align*}
u(1-u)\leq\frac{1}{4}, \qquad M\geq\frac{1}{2\gamma}.
\end{align*}
Combining \eqref{eq20} and \eqref{eq21}, we obtain
\begin{align}\label{eq22}
M'\geq\ep^{-\bar u}
\left(2M^2-\frac{M}{2\gamma}-\frac{c}{4\gamma^2}\right)
\quad\text{for}\quad M\geq\frac{1}{2\gamma}.
\end{align}
The expression in parentheses in \eqref{eq22} can be written as
\begin{align}\label{eq23}
2M^2-\frac{M}{2\gamma}-\frac{c}{4\gamma^2}
=2\left(M-\frac{1}{2\gamma}\right)^2
+\frac{3}{2\gamma}\left(M-\frac{1}{2\gamma}\right)
-\frac{c-1}{4\gamma^2}.
\end{align}
By \eqref{eq18}, we have
\begin{align}\label{eq24}
\frac{d}{dt}\left(\frac{\ep}{8\gamma}\log c(t)\right)
=-\frac{c(t)-1}{4\gamma^2}.
\end{align}
In view of \eqref{eq23} and \eqref{eq24}, let
\begin{align}\label{eq25}
\varphi(t)=M(t)-\frac{1}{\gamma}
\left(\frac{1}{2}+\frac{\ep}{8}\log c(t)\right).
\end{align}
For $\varphi>0$, we have 
\begin{align*}
M-\frac{1}{2\gamma}\geq\varphi>0.
\end{align*}
Combining \eqref{eq12} with \eqref{eq22}-\eqref{eq25}, we obtain
\begin{align}\label{eq26}
\varphi'
&=M'+\frac{c-1}{4\gamma^2}\nonumber\\
&\geq\ep^{-\bar u}
\left[2\left(M-\frac{1}{2\gamma}\right)^2
+\frac{3}{2\gamma}\left(M-\frac{1}{2\gamma}\right)\right]
+\frac{1-\ep^{-\bar u}}{4\gamma^2}(c-1)\nonumber\\
&\geq\frac{1}{\ep}\left(2\varphi^2+\frac{3}{2\gamma}\varphi\right)\quad\text{for a.e. }t\in[0,T). 
\end{align}

We claim that
\begin{align*}
\varphi(t)>0\quad\text{for all}\quad t\in[0,T).
\end{align*}
Otherwise, there exists $t_0<T$ such that
\begin{align*}
\varphi(t_0)=0,\qquad
\varphi(t)>0\quad\text{for all}\quad t\in[0,t_0).
\end{align*}
By \eqref{eq26}, we have
\begin{align*}
\varphi'(t)>0\quad\text{for a.e. }t\in[0,t_0).
\end{align*}
Hence
\begin{align*}
\varphi(t_0)\geq\varphi(0)=\delta_0>0.
\end{align*}
This contradicts $\varphi(t_0)=0$.
Consequently, we obtain
\begin{align*}
\varphi(t)\geq\varphi(0)=\delta_0>0
\quad\text{for all}\quad t\in[0,T),
\end{align*}
Hence, from \eqref{eq26} we have
\begin{align*}
\frac{d}{dt}\log\left(1+\frac{3}{4\gamma\varphi(t)}\right)
=-\frac{3\varphi'(t)}{\varphi(t)(3+4\gamma\varphi(t))}
\leq-\frac{3}{2\ep\gamma}
\quad\text{for a.e. }t\in[0,T).
\end{align*}
Integrating the above inequality from $0$ to $t$ yields
\begin{align*}
0<\log\left(1+\frac{3}{4\gamma\varphi(t)}\right)
\leq\log\left(1+\frac{3}{4\gamma\delta_0}\right)
-\frac{3t}{2\ep\gamma},
\end{align*}
which implies that the maximal existence time $T$ is finite
and satisfies \eqref{eq6}.
Moreover
\begin{align*}
M(t)\geq\varphi(t)\to+\infty
\quad\text{as }t\to T.
\end{align*}
This completes the proof of Theorem~\ref{thm:constant}.
\end{proof}

\subsection{Linear interaction potential}\label{subsec:linear}

\begin{proof}[\normalfont\bfseries Proof of Theorem~\ref{thm:linear}]
Along the characteristic flow \eqref{eq27}, let
\begin{align*}
g(t)=u_x(t,q(t,x_0)).
\end{align*}
Differentiating \eqref{eq3} with respect to $x$ yields
\begin{align}
\bar u_x
&=\frac{2}{\gamma^2}\int_x^{x+\gamma}u(y)\,dy-\frac{2}{\gamma}u(x),\label{eq28-1} \\ 
\bar u_{xx}
&=-\frac{2}{\gamma}u_x(x)
+\frac{2}{\gamma^2}\bigl(u(x+\gamma)-u(x)\bigr).
\label{eq28}
\end{align}
From \eqref{eq28-1} and $0\leq u\leq1$, we have
\begin{align}\label{eqbound}
-\frac{2u}{\gamma}\leq\bar u_x
\leq\frac{2(1-u)}{\gamma}.
\end{align}
Substituting \eqref{eq28} into \eqref{eq14}
and using $u(t,x+\gamma)\geq0$, we obtain
\begin{align}
g'(t)&\geq\ep^{-\bar u}\biggl[
2g^2+\left(2(1-2u)\bar u_x-\frac{2}{\gamma}u(1-u)\right)g-u(1-u)\left(\bar u_x^2+\frac{2u}{\gamma^2}\right)
\biggr]\nonumber\\
&=2\ep^{-\bar u}(g-g_1)(g-g_2),
\label{eq29}
\end{align}
where $g_2(\geq g_1)$ is given by
\begin{align*}
g_2=\frac12\biggl[
&-(1-2u)\bar u_x+\frac{u(1-u)}{\gamma}+\sqrt{
\left((1-2u)\bar u_x-\frac{u(1-u)}{\gamma}\right)^2
+2u(1-u)\left(\bar u_x^2+\frac{2u}{\gamma^2}\right)}
\biggr].
\end{align*}

We next estimate $g_2$. By \eqref{eqbound} and $u(1-u)\leq\frac{1}{4}$, we have
\begin{align}\label{eqq1-1}
(1-2u)\bar u_x-\frac{u(1-u)}{\gamma}
&\geq\frac{1}{\gamma}
\left((1-2u)^2-|1-2u|-\frac14\right)\notag\\
&=\frac{1}{\gamma}\left(|1-2u|-\frac12\right)^2
-\frac{1}{2\gamma}
\geq-\frac{1}{2\gamma},
\end{align}
and
\begin{align}\label{eqq2}
u(1-u)\left(\bar u_x^2+\frac{2u}{\gamma^2}\right)
&\leq\frac{1}{\gamma^2}
\left[4u(1-u)\bigl(u^2+(1-u)^2\bigr)+2u^2(1-u)\right]\notag\\
&\leq\frac{1}{\gamma^2}
\left[\frac12\bigl(2u(1-u)+u^2+(1-u)^2\bigr)^2+\frac{u}{2}\right]\notag\\
&=\frac{1+u}{2\gamma^2}
\leq\frac{1}{\gamma^2}.
\end{align}
Combining \eqref{eqq1-1} and \eqref{eqq2}, we have
\begin{align*}
g_2^2
&=\left(\frac{u(1-u)}{\gamma}-(1-2u)\bar u_x\right)g_2
+\frac{u(1-u)}2
\left(\bar u_x^2+\frac{2u}{\gamma^2}\right)\\
&\leq\frac{g_2}{2\gamma}+\frac1{2\gamma^2}.
\end{align*}
Solving this inequality yields
\begin{align}
g_2
\leq\frac12\left(
\frac1{2\gamma}
+\sqrt{\frac1{4\gamma^2}+\frac2{\gamma^2}}
\right)
=\frac1\gamma.
\label{eq30}
\end{align}

Since $g_1\leq0$, it follows from
\eqref{eq29}, \eqref{eq30} and \eqref{eq12} that
\begin{align}
g'(t)\geq\frac2{\ep}g(t)
\left(g(t)-\frac1\gamma\right)\quad\text{for}\quad g(t)\geq\frac1\gamma.
\label{eq31}
\end{align}

We claim that
\begin{align*}
g(t)>\frac1\gamma
\qquad\text{for all }t\in[0,T).
\end{align*}
Otherwise, there exists $t_0<T$ such that
\begin{align*}
g(t_0)=\frac1\gamma,\qquad
g(t)>\frac1\gamma\quad\text{for }t\in[0,t_0).
\end{align*}
Then \eqref{eq31} gives $g'>0$ on $[0,t_0)$,
and hence from condition \eqref{eq8}, we have
\begin{align*}
g(t_0)\geq g(0)>\frac1\gamma,
\end{align*}
which leads to a contradiction.
Consequently, we obtain
\begin{align*}
g(t)\geq g(0)>\frac1\gamma
\quad\text{for all }t\in[0,T).
\end{align*}
Integrating \eqref{eq31} from $0$ to $t$ yields
\begin{align*}
0<
\log\frac{\gamma g(t)}{\gamma g(t)-1}
\leq
\log\frac{\gamma u'_0(x_0)}{\gamma u'_0(x_0)-1}
-\frac{2t}{\ep\gamma},
\end{align*}
which implies that the maximal existence time $T$
is finite and satisfies \eqref{eq9}. Moreover,
\begin{align*}
u_x(t,q(t,x_0))=g(t)\to+\infty
\quad\text{as }t\to T.
\end{align*}
This completes the proof of Theorem~\ref{thm:linear}.
\end{proof}

\section*{Acknowledgment}
This work was supported by the National Natural Science Foundation of China (No. 12571261).

\end{document}